\documentclass[a4paper,10pt]{amsart}
\usepackage{graphicx,amssymb,amsmath,amsfonts,amsthm, color}
\usepackage{nicefrac}
\usepackage{geometry}
\usepackage{anysize}
\marginsize{2cm}{2cm}{2cm}{2cm}
\usepackage{mathtools}
\usepackage{textcomp}
\usepackage{paralist}
\usepackage[linkcolor=black,citecolor=black]{hyperref}
\usepackage[english]{babel}

\renewcommand {\epsilon}{\varepsilon}

\newcommand{\EE}{\mathbb{E}}

\newcommand{\NN}{\mathbb{N}}

\newcommand{\PP}{\mathbb{P}}

\newcommand{\RR}{\mathbb{R}}

\newcommand{\aA}{\mathcal{A}}
\newcommand{\bB}{\mathcal{B}}

\newcommand{\fF}{\mathcal{F}}

\newcommand{\jJ}{\mathcal{J}}
\newcommand{\kK}{\mathcal{K}}
\newcommand{\lL}{\mathcal{L}}

\newcommand{\oO}{\mathcal{O}}

\newcommand{\qQ}{\mathcal{Q}}

\newcommand{\sS}{\mathcal{S}}

\newcommand{\uU}{\mathcal{U}}

\newcommand{\wW}{\mathcal{W}}

\newcommand{\al}{\alpha}

\newcommand{\e}{\varepsilon}
\newcommand{\la}{\lambda}
\newcommand{\La}{\Lambda}

\newcommand{\om}{\omega}

\newcommand{\ra}{\rightarrow}

\newcommand{\ti}{\tilde}

\newcommand{\vzv}{\Leftrightarrow}

\newcommand{\zzv}{\Longrightarrow}

\newcommand{\lgl}{\ensuremath{\langle}}
\newcommand{\rgl}{\ensuremath{\rangle}}

\newcommand{\ind}{\mathbf{1}}

\newcommand{\lqq}{\leqslant}
\newcommand{\gqq}{\geqslant}

\newtheorem{thm}{Theorem}
\newtheorem{prop}{Proposition}
\newtheorem{cor}{Corollary}

\newtheorem{lem}{Lemma}

\newtheorem{rem}{Remark}
\newtheorem{exm}{Example}
\newtheorem{defn}{Definition}

\DeclareMathSymbol{\ophi}{\mathalpha}{letters}{"1E}

\renewcommand{\phi}{\varphi}

\newenvironment{proof2}{\par\noindent{\bf Proof\,:}}{\hfill$\blacksquare$\par}

\newfont{\cyrfnt}{wncyr10}
\def\J3{\cyrfnt{\rm \u{\cyrfnt I}}}
\def\j3{\cyrfnt{\rm \u{\cyrfnt i}}}

\usepackage[]{color}
\definecolor{DarkGreen}{rgb}{0.1,0.7,0.3}   

\definecolor{DarkGreen}{rgb}{0.1,0.7,0.3}   

\allowdisplaybreaks[4]

\begin{document}
\title[Moment estimates in the first Borel-Cantelli Lemma]{
  Moment estimates in the first Borel-Cantelli Lemma\\
with applications to mean deviation frequencies
}
\author{Luisa F. Estrada $\qquad$}
 \address{Departamento de Matem\'aticas, Universidad de los Andes, Bogot\'a,  Colombia. lf.estrada@uniandes.edu.co} 
\author{$\qquad$Michael A. H\"ogele}
\address{Departamento de Matem\'aticas, Universidad de los Andes, Bogot\'a, Colombia, ma.hoegele@uniandes.edu.co}
\hfill\\[-1cm]
\maketitle
    
\begin{abstract}
We quantify the elementary Borel-Cantelli Lemma by 
higher moments of the overlap count statistic in terms of the 
weighted summability of the probabilities. Applications include mean deviation frequencies in the Strong Law and the Law of the Iterated Logarithm. 
\end{abstract}
\hfill\\

\section{\textbf{Introduction}} 
    
\noindent The first Borel-Cantelli Lemma, which goes back to the seminal works of \'E. Borel \cite{Bo1909} and F. Cantelli \cite{Ca1917},
appears nowadays as one of the backbone results in the modern probability literature and is one of the essential tools for proving a.s. convergence results such as Etemadi's strong law of large numbers~\cite{Et81}. In its simplest and most widely used formulation it states that for a family of events $(E_n)_{n\in \NN}$ on a given probability space $(\Omega, \aA, \PP)$ the sufficiently fast decay of the probabilities $(\PP(E_n))_{n\in\NN}$, quantified by the finite summability 
\begin{equation}\label{first BC summability}C_1 = \sum_{n=1}^\infty \PP(E_n) < \infty,
\end{equation} 
yields that the overlap count $\oO := \sum_{n=1}^\infty\ind_{E_n}$ is a.s. finite. The famous one-line proof reads as follows: 
\begin{equation}\label{first BC}
\PP(\oO = \infty) = \PP(\limsup_{n\ra\infty} E_n) 
= \lim_{n\ra\infty} \PP(\bigcup_{m=n}^\infty E_m)  \lqq \lim_{n\ra\infty} \sum_{m= n}^\infty \PP(E_m) = 0.
\end{equation}
This result, however, is suboptimal in terms of moments since the monotone convergence theorem yields  
\begin{equation}\label{e:expectation}
\EE[\oO] = \EE[\lim_{n\ra\infty} \sum_{m=1}^n \ind_{E_m}] = \lim_{n\ra\infty} \sum_{m=1}^n \EE[\ind_{E_m}] = \sum_{m=1}^\infty \PP(E_m),
\end{equation}
whose right-hand side is finite by hypothesis. On the other hand, the law of $\oO$ has been explicitly known for a long time by the so-called Schuette-Nesbitt formula \cite{G79}, which gives a complete representation of the probabilities $\PP(\oO = k)$ for $k\in \NN_0$ in terms of a generalized inclusion-exclusion formula. Nevertheless, this formula includes the probabilities of all finitely many intersections of elements of $(E_n)_{n\in \NN}$, which is information that is rarely at hand in applications. Instead, the only information typically available, when applying the first Borel-Cantelli Lemma is the rate of convergence of the sequence $(\PP(E_n))_{n\in \NN}$, and in some cases, additionally, either the independence of the family  $(E_n)_{n\in \NN}$ or the monotonic inclusion $E_{n+1}\subset E_n$. 
This article addresses the following questions: 
\begin{enumerate}
 \item[i.)] Given \textit{only} the \textit{rate of convergence} $\PP(E_n)\to0$ as 
 $n\to\infty$, 
 what can be said about higher moments~of~$\oO$? 
 \item[ii.)] How can the results of i.) be improved by (a) the independence of the family  $(E_n)_{n\in \NN}$ or (b) the monotonicity $E_{n+1}\subset E_n$?   
\end{enumerate}

\noindent Roughly speaking, the body of literature on the Borel-Cantelli Lemmas \cite{Bo1909, Ca1917, Lev37, La63, Fr73, Ch12} concentrates along two lines of research: 1.)  extensions of the summability condition in different settings which imply $\oO< \infty$ a.s., or on some event of only positive probability such as in the Chung-Erd\H{o}s inequality~\cite{Amg06, CE51, DS67}, the Kochen-Stone lemma \cite{KS64} or the Erd\H{o}s-R\'enyi lemma \cite{ER59, Pe02, Pe06, Re70, FLS09}, and 2.) the precise asymptotic divergence of the overlap statistic in the case of nonsummable probabilities \cite{CE51, DF65, DS67,Br80, OW83, Hi83, Pe02, Amg06, Pe06, Xie09, FLS09, SXXZ09, XXXS13, Fr13}. 
To our knowledge, higher order moments of the overlap statistic of summable probabilities \eqref{first BC summability} have been studied in \cite{Fr73, Wi91} for martingale differences, which covers the independent case, but not that of arbitrary families of events. 
Very recently, we learned of the results in \cite{AO21+}, where the idea of studying the decay of $\sum_{m=n}^\infty\PP(E_m)$ combined with metastability appears, although not in the context of a systematic study of the moments of $\oO$. Moreover, recent applications of the overlap statistic, which include the Schutte-Nesbitt formula and moments of the overlap statistic in special cases, are found in information theory and finance \cite{KBK13, BK15, CDR13}. In the probability textbooks the Borel-Cantelli Lemmas appear essentially unchanged without taking advantage of the fine information in the Schuette-Nesbit formula \cite{G79}, see for example \cite{Bi99, Br68, Du10, Fe1-68, Ka02, Kle08, Sh99}. This short article offers a swift, elementary and useful improvement of the first Borel-Cantelli Lemma~\ref{e:expectation}, closing this somewhat surprising gap in the literature. Our main results include a simple quantification of higher order moments of the overlap statistic for general $(E_n)_{n\in \NN}$ (Theorem~\ref{prop: sufficient}) and we improve Freedman's universal upper bound (\cite{Fr73}, Prop. 17) in terms of the rate of convergence of the probabilities for independent $(E_n)_{n\in \NN}$ (Theorem~\ref{prop: main result 3}). In addition, we give a sample of the immediate usefulness of our findings in various statistical and theoretical applications. 

Due to the restriction on the rate of convergence in question i.), our results turn out to be slightly suboptimal, still, they provide a  quantitative improvement (in terms of the mean deviation frequency) of many standard results in the literature, that rely on the first Borel-Cantelli Lemma. For instance, we show a sharpened version of the \textit{second} Borel-Cantelli Lemma. We also introduce the notion of a.s. convergence with mean deviation frequency (MDF), which generalizes the concepts of complete convergence \cite{Yu99} and fast convergence \cite{Kle08}. This new concept is applied to different versions of the strong law of large numbers, such as the Vapnik-Chervonenkis (VC) theorem, the Glivenko-Cantelli theorem, Etemadi's strong law of large numbers, Cram\'er and Sanov's theorems with applications to the LDP of long rare sequences in \cite{DZ98} and the method of moments. Furthermore we refine the law of the iterated logarithm with an upper bound of the expected number of boundary crossings and an application to MDF errors for strong numerical schemes of SDE of higher order. All the results we have mentioned allow for rather straightforward goodness-of-fit test procedures.

\section{\textbf{The main results: higher order moments in the first Borel-Cantelli Lemma}} 

\noindent 
This section provides the moment estimates of the overlap statistic in the setting of the first Borel-Cantelli Lemma in Subsection~\ref{2.1}, followed by the setting of the second Borel-Cantelli Lemma in Subsection~\ref{2.2}. 

\subsection{Higher order moments for the overlap statistic in the first Borel-Cantelli Lemma}\label{2.1}\hfill\\[-3mm]

\noindent For the setting ii.(b) of a monotonic family of events $E_n \supset E_{n+1}$ we have $\{\oO = k\} = E_k \setminus E_{k+1}$ such that the moments can be calculated expicitly. Due to the monotonicity there are no gaps, which implies that the values of $\oO$ turn out to be maximal values to which general families of events can be compared to in Theorem~\ref{prop: sufficient}. 

\subsubsection{\textnormal{\textbf{Higher order moments for the overlap statistic of nested events}}}
\begin{prop}\label{prop: nested}
Let $(\Omega, \aA, \PP)$ be a probability space and $(E_n)_{n\in \NN_0}$ a nested family of events, $E_n\supset E_{n+1}$, $n\in \NN_0$.
Then for any given nonnegative sequence $(a_n)_{n\in \NN_0}$ with $\sS(N) := \sum_{n=0}^N a_n$, $N\in\NN$, we have that 
\begin{align}\label{e:genmoments0}
\EE[\sS(\oO)] = \sum_{n=0}^\infty a_n \cdot\PP(E_n),\qquad \mbox{ whenever the right-hand side is finite.}  
\end{align}
\end{prop}

\begin{proof} 
    By construction we have $\{\oO = k\} = E_k \setminus E_{k+1}$, $k\in \NN_0$. Fix $N\in \NN$ and $\oO_N := \sum_{n=0}^N \ind_{E_n}$. Then we have $\{\oO = k\} = E_k \setminus E_{k+1}$ for $k = 0,\dots, N-1$, while $\{\oO_N = N\} = E_N$. Note that $\{\oO = 0\} = \{\oO_N = 0\} = \Omega \setminus E_0$. 
    Let us denote $p_k = \PP(E_k)$, $k\in \NN_0$, and formally $p_0 =1-\PP(E_0)$. 
    We apply the elementary summation by parts formula for sequences $(f_k)_{k\in \NN}$ and $(g_k)_{k\in \NN}$  
    \begin{align}\label{e: summationbyparts}
    \sum_{n=0}^N f_n g_n - f_N \sum_{n=0}^N g_n = - \sum_{k=0}^{N-1} (f_{k+1}-f_k) \sum_{n=0}^k g_n, \qquad N\in \NN. 
    \end{align}
    For $f_n = p_n$ and $g_n = a_n$, and due to $\PP(\oO_N = 0) = 1- p_0$, $\PP(\oO_N = k) = p_k - p_{k+1}$, $\PP(\oO_N = N) = p_N$ we have 
    \begin{align}\label{e:sumbyparts0}
    \sum_{n=0}^N a_n \,p_n -p_N \cdot\sS(N) = \sum_{k=0}^{N-1} (p_{k}-p_{k+1}) \cdot \sS(k) = \sum_{k=0}^{N-1} \sS(k)\cdot \PP(\oO_N = k)= \EE[\sS(\oO_N)]-p_N \cdot\sS(N), N\in \NN_0. 
    \end{align}
    We cancel $p_N \cdot\sS(N)$ and conclude by sending $N\ra\infty$ with the help 
    of the monotone convergence theorem. 
\end{proof}

\subsubsection{\textnormal{\textbf{The first main result: Higher order moments of the overlap of general families of events}}} We now treat question i.) and apply Proposition~\ref{prop: nested} to families of non-nested events. 

\begin{thm}[Higher order moments of the overlap statistic $\oO$]\label{prop: sufficient}
Let $(\Omega, \aA, \PP)$ be a probability space and $(E_n)_{n\in \NN_0}$ an arbitrary  family of events. 
Then for any given nonnegative sequence $(a_n)_{n\in \NN_0}$ with $\sS(N) := \sum_{n=0}^N a_n$, $N\in \NN$, we have that 
\begin{align}\label{e:genmoments}
\EE[\sS(\oO)]\lqq \sum_{n=0}^\infty a_n \sum_{m=n}^\infty \PP(E_m), \qquad 
\mbox{ whenever the right-hand side is finite. }
\end{align}
    \end{thm}

\begin{proof2}
Define $\ti E_n := \bigcup_{m=n}^\infty E_n$, $n\in\NN_0$. By construction, $E_n\subset \ti E_n$, $n\in \NN_0$, and the events $(\ti E_n)_{n\in \NN_0}$ are nested in the sense of $\ti E_{n}\supset \ti E_{n+1}$, $n\in \NN_0$. The monotonicity $E_n\subset \ti E_n$ implies that $\ind_{E_n} \lqq \ind_{\ti E_N}$ a.s. and 
\[
\oO =\sum_{n=0}^\infty \ind_{E_n}\lqq \sum_{n=0}^\infty \ind_{\ti E_n} =: \ti \oO\quad \mbox{a.s.}
\]
Since $(\sS(N))_{N\in \NN_0}$ is nondecreasing, Proposition~\ref{prop: nested} combined with $\PP(\ti E_n) \lqq \sum_{m=n}^\infty \PP(E_m)$, $n\in \NN_0$, yields
\begin{align*}
\EE[\sS(\oO)] \lqq \EE[\sS(\ti \oO)] = 
\sum_{n=0}^\infty a_n \cdot\PP(\ti E_n)\lqq \sum_{n=0}^\infty a_n \cdot \sum_{m=n}^\infty \PP(E_m).\\[-1.2cm] 
\end{align*}
\end{proof2}

\begin{cor}[Polynomial and exponential moments of the overlap statistic $\oO$]\label{cor: sufficient poly and exp}
Let $(\Omega, \aA, \PP)$ be a probability space and $(E_n)_{n\in \NN}$ an arbitrary family of events $(E_n)_{n\in \NN}$. Then we have: 
\begin{align*}
&\mbox{\textnormal{(1)}} \quad K_1(p):= \sum_{n=1}^\infty n^p \sum_{m=n}^\infty \PP(E_m) <\infty\quad \mbox{ for }p>0  \qquad \zzv \qquad \EE[\oO^{p+1}] \lqq (p+1)\, K_1(p).\\
&\mbox{\textnormal{(2)}} \quad K_2(p):= \sum_{n=0}^\infty e^{np} \sum_{m=n}^\infty \PP(E_m) <\infty \quad \mbox{ for }p>0  \qquad \zzv \qquad \EE[e^{p\oO}] \lqq K_2(p) +1.\\[-0.5cm] 
\end{align*}
\end{cor}

\begin{proof2} For simplicity we show the results for $p\in\NN$. The general case follows by interpolation. 
For $a_n = n^p$, $p\in \NN$, the sum of monomials satisfies 
\begin{equation}\label{e:Faulhaber}
\sS(N) = \sum_{n=1}^N n^p = \frac{(N+1)^{p+1}}{p+1} + \frac{1}{2} N^p + \sum_{j=2}^p \binom{p}{j} \frac{B_j}{p+1-j} N^{p+1-j} =: \fF_{p+1}(N), \quad N\in\NN, 
\end{equation} 
where $\fF_{p+1}(N)$ is the Faulhaber polynomial of order $p+1$ 
and $B_{j}, j=2, \dots, p$ are the Betti-numbers (see \cite{JB1713, Jac1834}). 
Then Theorem~\ref{prop: sufficient} implies 
$\EE[\fF_{p+1}(\oO)] \lqq K_1(p)$, where $K_1(p)$ is finite by hypothesis. Moreover, the leading coefficent of the Faulhaber polynomials is known to be $1 / (p+1)$ and,  since all coefficients of $\fF_{p+1}$ are positive, by monotonicity we obtain that 
\[
\EE[\oO^{p+1}] \lqq (p+1) \EE[\fF_{p+1}(\oO)] \lqq (p+1) K_1(p),  
\]
which proves item (1). For $a_n = e^{np}$ we have the geometric sum 
$\sS(N) := \sum_{n=0}^N a_n = (e^{(N+1)p}-1)/(e^p-1).$
Consequently, Theorem~\ref{prop: sufficient} yields 
\begin{align*}
\EE\Big[\frac{e^{p(\oO+1)}-1}{e^p-1}\Big] \lqq K_2(p), \mbox{ and item (2) follows from } \EE[e^{p\oO}] 
= \frac{(e^p-1)}{e^p} \EE\Big[\frac{e^{p(\oO+1)}-1}{e^p-1}\Big]+e^{-p}\lqq K_2(p)+1.
\end{align*}
\end{proof2}

\begin{exm}\label{ex: poly}
Assume $\PP(E_n) \lqq C / n^q$, $n\in \NN$, for some constants $C, q>0$. 
Then we have for all $p+1 < q-1$ 
\begin{equation}\label{e: Riemann}
\EE[\oO^{p+1}] \lqq C \cdot \zeta(q-1-p), \qquad \mbox{ where }\zeta(z) := \sum_{n=1}^\infty \frac{1}{n^z} \mbox{ is Riemann's zeta function.} 
\end{equation} 
The direct calculation of the constant $K_1(p)$ in Corollary~\ref{cor: sufficient poly and exp}, item~(1), yields that for $p < q-2$ 
we have 
\[
K_1(p) = \sum_{n=1}^\infty n^{p} \sum_{m=n}^\infty \PP(E_m)
\lqq C \sum_{n=1}^\infty n^{p} \int_n^\infty \frac{1}{x^{q}} dx =\frac{C}{q-1} 
\zeta(q-1-p)<\infty,  
\]
which consequently implies \eqref{e: Riemann}. Note that the condition $q>p+2$ is clearly suboptimal, since for any $q>1$ we already know by \eqref{e:expectation} that moments of order $p=1$ are finite. 
\end{exm}

\begin{exm}\label{exm: exp}
Assume $\PP(E_n) \lqq C b^n$ for all $n\in \NN_0$ for some $b\in (0, 1)$ and $C>0$.  Then Corollary~\ref{cor: sufficient poly and exp}, item~(2), implies for all $0 < p< |\ln(b)|$ the estimate 
\begin{equation}\label{e:expmoment}
\EE[e^{p\oO}] \lqq 1+ \sum_{n=0}^\infty e^{pn} 
\sum_{m=n}^\infty \PP(E_m) \lqq 1+ 
\sum_{n=0}^\infty e^{pn} \frac{C b^n}{1-b} = 1+ \frac{C}{1-b} \sum_{n=0}^\infty e^{n (p + \ln(b))} = \frac{C}{(1-b)(1- e^p b)}  +1.
\end{equation}
\end{exm}

\subsection{Exponential moments of the overlap in the counterpart of the second Borel-Cantelli Lemma}\label{2.2}\hfill\\[-3mm]

\noindent The following classical result relies on Kolmogorov's three series theorem (see e.g.   \cite{Wi91},~Sec.~12.4). 
 
\begin{lem}\label{lem:independent} 
Let $(\Omega, \aA, \PP)$ be a probability space and $(E_n)_{n\in \NN}$ an independent family of events. 
Then  
\begin{align*}
\oO < \infty \mbox{ a.s.} \qquad \mbox{ is equivalent to } \qquad \mbox{Var}(\oO) < \infty \qquad \mbox{ and also equivalent to }\qquad \eqref{first BC summability}. 
\end{align*}
In particular, any of the preceding conditions implies $\EE[\oO^2] \lqq  C_1 \cdot(1 + C_1)<\infty$. 
\end{lem}

\noindent Lemma~\ref{lem:independent}, however, is suboptimal in terms of integrability of $\oO$. Freedman showed in \cite{Fr73} that under the same conditions all exponential moments are finite with explicitly known optimal universal upper bounds. 

\subsubsection{\textnormal{\textbf{Freedman's universal bound of the exponential moments in the second B.-C. Lemma}}} In Proposition 17 of \cite{Fr73} Freedman shows the following result in a more general martingale differences setting. We give an elementary alternative proof. 

\begin{thm}[Freedman's universal bound, \cite{Fr73}]
\label{prop: universalexp}
Let $(\Omega, \aA, \PP)$ be a probability space and $(E_n)_{n\in \NN}$ an independent family of events. 
Then the Borel-Cantelli summability \eqref{first BC summability} with constant $C_1>0$ implies 
\[
\EE\big[e^{r \oO}\big] \lqq \exp\big(C_1 (e^r-1)\big)\qquad \mbox{ for all }r>0.  
\]
In particular, for all $k\in \NN$ we have 
\begin{equation}\label{e: FreedProb}
\PP(\oO\gqq k)\lqq \inf_{r>0} \exp(-k r+ C_1(e^r -1)) = \exp(- k\ln(k) + k (\ln(C_1)+1) - C_1).
\end{equation}
\end{thm}

\begin{proof2} The Bernoulli inequality $1+x \lqq e^x$, $x\in \RR$, implies for any $r>0$  
\begin{align*}
\EE\big[e^{r \oO}\big]& = \prod_{n=1}^\infty \EE[e^{r \ind_{E_n}}] = \prod_{n=1}^\infty \big(e^r \PP(E_n) + 1- \PP(E_n)\big)  
= \prod_{n=1}^\infty \exp\big(\ln(1+ e^r \PP(E_n) - \PP(E_n))\big) \\ 
&= \exp\big(\sum_{n=1}^\infty \ln(1+ (e^r  - 1) \PP(E_n))\big)  
\lqq \exp\big( (e^r-1) \sum_{n=1}^\infty \PP(E_n)\big)  = \exp\big(C_1(e^r-1)\big).
\end{align*}
Markov's inequality with a subsequent minimization yields \eqref{e: FreedProb}. 
\end{proof2}
\noindent 
As a consequence of Theorem \ref{prop: universalexp} we formulate a version of the second Borel-Cantelli Lemma which illustrates the sharp dichotomy between summability and nonsummability of $\sum_{n=1}^\infty \PP(E_n)$ in terms of the moments of $\oO$.

\begin{cor}[Freedman's universal moment version of the second Borel-Cantelli Lemma] 
Let $(\Omega, \aA, \PP)$ be a probability space and $(E_n)_{n\in \NN}$ an independent family of events. 
Then we have: 
\begin{align}
\mbox{1. }&\sum_{n=1}^\infty \PP(E_n) = \infty \quad \zzv\quad \oO = \infty \quad \mbox{a.s.}\nonumber\\[-1mm]
\mbox{2. } &\sum_{n=1}^\infty \PP(E_n) <\infty \quad \zzv\quad \oO < \infty \quad \mbox{a.s. with } ~~\EE[e^{r\oO}] \lqq \exp\Big((e^r-1)\cdot \sum_{n=1}^\infty \PP(E_n)\Big) \lqq \infty ~\mbox{ for all }r>0.\label{e: universal}
\end{align}
\end{cor}
\noindent 
Note that the upper bound \eqref{e: universal} does not depend on the rate of convergence of $(\PP(E_n))_{n\in \NN}$, but only on the value of the sum $C_1$ in \eqref{first BC summability}. In Corollary~\ref{cor: specificexpmoments} of Subsection~\ref{ss:improved} we improve Freedman's result to smaller upper bounds which depend on the rate of convergence of 
$(\PP(E_n))_{n\in\NN}$.

\subsubsection{\textnormal{\textbf{The second main result: Improved exponential moment bounds in the second BC-Lemma}}}\label{ss:improved} We now address question ii.(a). The Schuette-Nesbitt formula \cite{G79} 
provides upper bounds for the exponential moments of $\oO$ using the rate of convergence of the tail function $\NN\ni m\mapsto L(m) = \sum_{n=m}^\infty\PP(E_n)$.

 \begin{thm}\label{prop: main result 3}
 Given $(\Omega, \aA, \PP)$ and $(E_n)_{n\in \NN}$ an independent family of events satisfying \eqref{first BC summability} with $C_1 < 1$. 
Then  
\[\EE\big[e^{r \oO}\big] \lqq (1-C_1 e^{r})^{-1}\qquad \mbox{ for all }r <|\ln(C_1)|.\]  
\end{thm}

 \begin{proof2}
For $0\lqq k\lqq N$, $k, N\in \NN$,  consider $G_k^N := \{\oO_N = k\}$ and $\oO_N := \sum_{n=1}^N \ind_{E_n}$. By the Schuette-Nesbitt formula~\cite{G79}, we have for any nonnegative sequence $(a_n)_{n\in \NN}$ the equality  
\begin{align}
\sum_{k=0}^N a_k\,\PP(G_k^N)  = \sum_{n=0}^N \qQ_n^N (a_n - a_0),\qquad \mbox{ where }\quad \qQ_n^N =  \sum_{\substack{J\subset \{1, \dots , N\}\\ | J | = n }} \PP(\bigcap_{j\in J} E_j).\label{e: SN}
\end{align}
Using the independence of $(E_j)_{j\in \NN}$ we get for all $n\in \NN$, $n\lqq N$,  
\begin{align*}
\qQ_n^N &=  \sum_{\substack{J\subset \{1, \dots , N\}\\ | J | = n }} \prod_{j\in J} \PP(E_j)
= \sum_{i_1=1}^{N-n} \sum_{i_2=i_1 +1}^{N-n+1} \hdots \sum_{i_{n-1} = i_{n-2} +1}^{N-1} \sum_{i_{n} = i_{n-1} +1}^{N} \prod_{\ell= 1}^n \PP(E_{i_\ell}) \lqq \Big(\sum_{i=1}^N \PP(E_i)\Big)^n\lqq C_1^n. 
\end{align*}
Hence \eqref{e: SN} yields for $a_k = e^{rk}$ 
with $r <|\ln(C_1)|$ for all $N\in \NN$ 
\begin{align}
\sum_{k=0}^\infty \PP(G_k^N) \cdot e^{rk} = \sum_{k=0}^N \PP(G_k^N) \cdot e^{rk} = \sum_{n=0}^N \qQ_n^N (e^{rn} - 1) \lqq  
\sum_{n=0}^N e^{r n} \qQ^N_n 
&\lqq \sum_{n=0}^N e^{rn}  C_1^n \lqq (1- C_1 e^r)^{-1}.\label{e:prefatou}
\end{align}
In addition, $G^N_k = G^{N+1}_{k} \cup G^{N+1}_{k+1}$ yields that $\PP(G^N_k)\gqq \PP(G^{N+1}_k)\gqq 0$ for all $N\gqq k\gqq 0$ such that the limit $\lim_{N\ra\infty} \PP(G^N_k)$ exists and is equal to $\PP(G_k)$ by measure continuity from above. Combining the preceding limit with \eqref{e:prefatou} and Fatou's Lemma finishes the proof, since
\begin{align*}
\sum_{k=0}^\infty \PP(G_k) \cdot e^{rk} = 
\sum_{k=0}^\infty \liminf_{N\ra\infty} \PP(G^N_k) \cdot e^{rk} \lqq 
\liminf_{N\ra\infty} \sum_{k=0}^\infty \PP(G^N_k) \cdot e^{rk} 
\lqq (1- C_1 e^r)^{-1}.\\[-1cm]
\end{align*}
 \end{proof2}

\begin{cor}\label{cor: specificexpmoments}
For a probability space $(\Omega, \aA, \PP)$ and an independent family of events $(E_n)_{n\in \NN}$ satisfying~\eqref{first BC summability} let  
\begin{equation}
C_m := \sum_{n=m}^\infty \PP(E_n), \quad m\in \NN, \qquad \mbox{ and }\qquad 
N_r(\delta) := \inf\{m\in \NN~|~ C_m < e^{-r}/\delta\}, \quad \mbox{ for }r>0, ~\delta>1.
\end{equation}
Then for all $r>0$ we have  
\begin{align*}
\EE\big[e^{r \oO}\big] \lqq \inf_{\delta>1}\inf_{m\gqq N_r(\delta)} e^{r m} (1- C_m e^r)^{-1} = \inf_{\delta>1}\frac{\delta}{\delta-1} e^{r \cdot N_r(\delta)} 
\end{align*}
and for any nonincreasing, invertible function $L: (0, \infty) \ra (0,\infty)$ such that $L(m) = C_m$, $m\in \NN$, we have 
\[
\EE\big[e^{r \oO}\big] \lqq  \inf_{\delta>1}\frac{\delta}{\delta-1}  \exp\big(r \cdot L^{-1}(e^{-r}/\delta)\big)\qquad \mbox{ for all } r>0. 
\]
\end{cor}
\noindent The proof relies on the fact that  $N_r(\delta)$ is the minimal $m\in \NN$ such that   
$\delta C_{m} <e^{-r}$ combined with Theorem~\ref{prop: main result 3}.

\begin{cor}[Improved exponential moment bounds in the second Borel-Cantelli Lemma] 
Let $(\Omega, \aA, \PP)$ be a probability space and $(E_n)_{n\in \NN}$ an independent family of events. 
Then we have: 
\begin{align*}
\mbox{1. }&\sum_{n=1}^\infty \PP(E_n) = \infty \quad \zzv\quad \oO = \infty \quad \mbox{a.s.}\\[-1mm]
\mbox{2. } &\sum_{n=1}^\infty \PP(E_n) <\infty \quad \zzv\quad \oO < \infty \quad \mbox{a.s.}\quad \mbox{ with } \quad \EE\big[e^{r \oO}\big] \lqq \inf_{\delta>1}\frac{\delta}{\delta-1}\, \exp(r \cdot L^{-1}(e^{-r}/\delta)) < \infty
\\[-3mm]
&\hspace{-2mm} \mbox{ for all }r>0\mbox{ and any } L: (0, \infty) \ra (0, \infty) \mbox{ nonincreasing, invertible such that }  L(m) = \sum_{n=m}^\infty \PP(E_n), m\in \NN. \end{align*}
\end{cor}
\noindent The proof is a direct application of Corollary~\ref{cor: specificexpmoments}. 

\begin{exm}\label{ex: indexpmoments}
For an independent family $(E_n)_{n \in\NN}$ with $\PP(E_n) \lqq c / n^p$, $n\in \NN$, and some $c>0$, $p>1$,  consider 
 $L(r) = c/ r^p$, $L^{-1}(s) = \big(c / s\big)^{1/p}$ and $L^{-1}(e^{-r}/\delta) = (\delta c\, e^{r})^{1/p}$ for $s,r>0$, $\delta>1$. Corollary~\ref{cor: specificexpmoments} yields 
\begin{equation}\label{e: Tupperbound}
\EE\big[e^{r \oO}\big] \lqq  \inf_{\delta>1}\frac{\delta}{\delta-1} \exp\big((\delta c)^{1/p}\cdot r\, e^{r/p}\big)\lqq 2 \exp\big((2 c)^{1/p}\cdot r \, e^{r/p} \big)\qquad \mbox{ for all } r>0. 
\end{equation}
Note that for values $p, \delta>1$ close to $1$ the exponential order of $r$ on the right-hand side of \eqref{e: Tupperbound} approaches the universal upper bound  given in 
Theorem~\ref{prop: universalexp}, i.e. the universal upper bound \eqref{e: universal} is optimal. In fact, the optimality statement is proven rigorously
in Proposition~(17) of Freedman~\cite{Fr73}. 
By Markov's inequality we have 
\[
\PP(\oO\gqq k)\lqq 
\inf_{r>0} 2 \exp\big(-k r + (2 c)^{1/p} r e^{r/p}\big) 
\lqq 2 \exp\big(\inf_{r>0} (-k r + (2 c)^{1/p} r e^{r/p})\big), \qquad k\in \NN.
\]
The minimizer is $r=p(\wW(k / (e(\delta c)^{1/p}) ) - 1)$, 
where $\wW$ is the principal branch of Lambert's $\wW$ function with the well-known asymptotics $\wW(x) = \ln(x) - \ln(\ln(x)) + o(1)_{x\ra\infty}$ (see \cite{CGHJK96}). This implies for $k\in \NN$, $k> e^2$, that 
\begin{align*}
\inf_{r>0}(-k r + (2 c)^{1/p} r e^{r/p}) &=  -pk \frac{(\wW(k / (e(2 c)^{1/p}) ) - 1)^2}{\wW(k/(e(2 c)^{1/p}))} 
~ =   
-p k \ln(k) + pk \ln(\ln(k)) + o(1)_{k\ra\infty}. 
\end{align*}
Hence, there is a constant $\kK = \kK(\delta, p, c)>0$ such that 
$
\PP(\oO\gqq k)\lqq \kK\cdot\exp(-p k [\ln(k) - \ln(\ln(k))])$, $k>e^2$.
\end{exm}

\begin{exm}
For an independent family $(E_n)_{n \in\NN}$ with $\PP(E_n) \lqq c\cdot b^n$, $n\in \NN$, 
for some constants $0< b < 1$ and $c>0$ we calculate for 
$L(r) = c b^r, r>0,$ its inverse $L^{-1}(s) = \ln( s/c) / \ln(b), s>0,$ such that for $\delta>1$  
\[L^{-1}(e^{-r}/\delta) 
=  (r +\ln(\delta c))/|\ln(b)|.\]
Therefore, Corollary~\ref{cor: specificexpmoments} implies for all $r>0$ 
\begin{equation}\label{ex: indepupperbound}
\EE\big[e^{r \oO}\big] \lqq  \inf_{\delta>1}\frac{\delta}{\delta-1} \exp\big(\big[r^2  + r \cdot \ln(\delta c)\big] / |\ln(b)|\big)\lqq 2 
\exp\big(\big[r^2  + r \cdot \ln(2 c)\big] / |\ln(b)|\big).
\end{equation}
By Markov's inequality and a subsequent minimization procedure we obtain the Gaussian type decay for $k\in \NN$
\[
\PP(\oO\gqq k) \lqq 2 \inf_{r>0} \exp\big( \big(r^2 + r \cdot [\ln(2c) -k |\ln(b)|\big]\big)/|\ln(b)|\big)
=2\exp\big(-(|\ln(b)| /4) \big[k  - (\ln(2C)/|\ln(b)|)\big]^2 \big). 
\]
This upper bound is much smaller than the asymptotics obtained in Example~\ref{ex: indexpmoments}.
\end{exm}

\section{\textbf{Applications to mean deviation frequencies (MDF)}} 

\subsection{Almost sure convergence with higher order MDF}\hfill\\[-3mm]

\noindent In general, the notion of a.s. convergence has the disadvantage that it is not easily grasped statistically, since the modulus of convergence and the number of deviations are random. In the light of Theorems~\ref{prop: sufficient}, ~\ref{prop: universalexp} and ~\ref{prop: main result 3} it is natural to introduce the notion of a.s. convergence distinguishing different moments of the deviation frequency.

\begin{defn}\label{def1}
Given a probability space $(\Omega, \aA, \PP)$, a sequence of random vectors $(X_n)_{n\in \NN}$, a random vector $X$ and $\e>0$ consider a nondecreasing function $\Lambda_\e: \NN\ra [0, \infty)$.  
\[
\mbox{\normalfont If (I)}\lim\limits_{n\ra\infty} X_n = X \mbox{ a.s.}\quad \mbox{ and }\quad 
\mbox{\normalfont(II)} \quad \EE[\Lambda_\e(\oO_\e)]< \infty, \quad \mbox{ for }\quad  \oO_\e := \sum_{n=1}^\infty \ind{\{|X_n-X|>\e\}}, \quad \e\in (0, 1). 
\]
we say that \textbf{$(X_n)_{n\in \NN}$ converges to $X$ 
a.s. with mean deviation frequency (MDF) of order $\La_\e$}.  
\end{defn}

\noindent A classical step to infer a.s. convergence is to show convergence in probability and subsequently to strengthen it by the summability of the error estimates and the first Borel-Cantelli Lemma \eqref{first BC}. Theorem~\ref{prop: sufficient} allows to quantify the excess of bare summability in terms of the moments of the number of deviations. 
Yukich introduced the notion of complete convergence in~\cite{Yu99}, which coincides with the a.s. convergence with MDF of order $1$. 
    
\begin{cor}[Complete convergence implies first order MDF]\label{cor:complete}
Given a probability space $(\Omega, \aA, \PP)$, a sequence of random vectors $(X_n)_{n\in \NN}$ and a random vector $X$. Then the complete convergence of $X_n\ra X$ is given by 
\begin{equation}\label{e: summability}
\phi(\e) := \sum_{n=1}^\infty  \PP(|X_n-X|\gqq \e) <\infty,\quad \e\in (0, 1), \quad \mbox{ and implies }\quad \EE[\oO_\e] = \phi(\e). 
\end{equation}
In particular, by the Markov inequality we have that $\PP(\oO_\e \gqq k) \lqq \phi(\e)/k$ for all $k\in \NN.$
\end{cor}
\noindent Corollary~\ref{cor:complete} is a direct consequence 
of the expectation shown in \eqref{e:expectation}. 

\begin{rem}[Fast convergence implies first order MDF]\label{rem: fast}
Due to the Markov inequality (see Theorem 6.12 (i) \cite{Kle08}) a sufficient condition for \eqref{e: summability} is $\sum_{n=1}^\infty  \EE[|X_n-X|^p] <\infty.$ 
\end{rem}

\noindent The consistency results of Corollary~\ref{cor:complete} and Remark~\ref{rem: fast}, allows us to quantify higher MDF. 

\begin{cor}[a.s. convergence with polynomial MDF]\label{cor: aspoloutlier} 
Given a probability space $(\Omega, \aA, \PP)$, a sequence of random vectors $(X_n)_{n\in \NN}$ and a random vector $X$. We assume the existence of a nondecreasing function $p: (0, 1) \ra [0, \infty)$ 
such that for $\La_\e(n) = n^{p(\e)+1}, n\in \NN, \e\in (0,1)$ we have  
\[
\phi(\e) := \sum_{n=1}^\infty n^{p(\e)} \sum_{m=n}^\infty \PP(|X_m-X|\gqq \e) <\infty, \qquad \mbox{ for all }\e\in (0, 1),   
\]
Then Corollary~\ref{cor: sufficient poly and exp} implies that $X_n \ra X$ a.s. with polynomial MDF of order $\La_\e$, i.e. for all $\e\in (0, 1)$ it follows  
\[
\EE[\oO_\e^{p(\e)+1}] \lqq (p(\e)+1) \cdot \phi(\e), 
\quad \mbox{and in particular }\quad 
\PP(\oO_\e\gqq k) \lqq k^{-(p(\e)+1)}\cdot(p(\e)+1)  \phi(\e), \qquad k \in\NN.
\]
\end{cor}
\noindent The proof is a direct application of Corollary~\ref{cor: sufficient poly and exp}, item~(1) and Definition~\ref{def1}.

\begin{cor}[a.s. convergence with exponential MDF]\label{cor: asexpoutlier}
Given a probability space $(\Omega, \aA, \PP)$, a sequence of random vectors $(X_n)_{n\in \NN}$ and a random vector $X$. If $X_n \ra X$ in probability and there exists a nondecreasing function $p: (0, 1) \to (0, \infty)$ such that for $\La_\e(n) :=\exp(n p(\e))$, $n\in \NN$, $\e\in (0,1)$, we have that 
\[
\phi(\e) := \sum_{n=0}^\infty \exp\big(n p(\e)\big) \sum_{m=n}^\infty \PP(|X_m-X|\gqq\e) <\infty, \qquad \e\in (0, 1). 
\]
Then Corollary~\ref{cor: sufficient poly and exp} implies that $X_n \ra X$ a.s. with exponential MDF of order $\La_\e$. That is, for all $\e \in (0,1)$ we have $\EE\big[e^{p(\e) \oO_\e}\big] \lqq \phi(\e)+1$, 
and $\PP(\oO_\e\gqq k) \lqq \exp(-p(\e) k) (\phi(\e)+1)$ for all $k\in\NN_0$. 
\end{cor}
\noindent The proof is a direct application of Corollary~\ref{cor: sufficient poly and exp}, item~(2) and Definition~\ref{def1}.

\medskip
\subsection{Mean deviation frequency in the Strong Law and Large Deviations}\hfill\\[-2mm]

\noindent The direct proof of the strong law of large numbers (SLLN) is based on a strengthening the weak law of large numbers by the Borel-Cantelli Lemma, see \cite{Et81}. Kolmogorov's strong law characterizes the validity of the SLLN in the case of i.i.d. sample means by the existence of the first moments. However, if the summands have higher order moments a stronger summability of the error probabilities in the weak law of large numbers is attained. In this setting, Theorem~\ref{prop: main result 3} yields estimates on the moments of the deviation frequency of a.s. convergence. This result is naturally improved in the presence of a large deviations principle. 

\subsubsection{\textnormal{\textbf{Mean deviation frequency estimate 
for relative frequencies}}}

In this section we quantify the a.s. convergence in the classical VC theorem and the Glivenko-Cantelli theorem.  We cite Theorem 2 in~\cite{VC71}. 

\begin{thm}[Classical VC inequality]\label{thm:VCclassical} 
Consider a VC class $S$ of sets in $\RR^d$ 
with index $m^S$, see \cite{VC71}.
The probability of the event $\pi^{\ell}$ 
that the relative frequency of at least one event in class $S$ 
differs from its probability in an experiment of size by more 
then $\e>0$, for $\ell \gqq 2/\e^2$, satisfies 
\begin{equation}\label{e: VC}
\PP(\pi^{(\ell)} > \e) \lqq 4 m^S(2\ell) e^{-\e^2 \ell /8}, \qquad \ell\in \NN.  
\end{equation}
\end{thm}
\noindent Combining \eqref{e: VC} with the classical Borel-Cantelli lemma 
the authors prove Theorem~3 in \cite{VC71}, which states that $m^S(\ell)\lqq \ell^p +1$, $\ell\in \NN$, for some $p>0$ implies that $\pi{(\ell)} \ra 0$ $\PP$-a.s.
For $\e>0$ denote $\oO_\e := \sum_{\ell=1}^\infty \ind\{\pi^{(\ell)}>\e\}$. 
\begin{thm}[VC theorem with Gamma type MDF] 
Under the hypotheses of Theorem \ref{thm:VCclassical} assume that $m^S(\ell)\lqq \ell^p +1$, $\ell\in \NN$, for some $p>0$. 
Then for any $\e, \delta>0$ and $N\gqq 2/\e^2$ we have for 
\[
\Lambda_{\e, \delta}(N):=  \frac{e^{\e^2 N /8}}{N^{1+\delta} m^S(2N)} \qquad \mbox{ that  }\qquad \pi{(\ell)} \ra 0\qquad \mbox{a.s.  with exponential MDF }\quad \EE[\Lambda_{\e, \delta}(\oO_\e)] <\infty. 
\]
In particular, we have the Gamma-function like asymptotics 
$\PP(\oO_\e\gqq k)\lqq C / \La_{\e, \delta}(k)$, $k\in \NN$ for some $C>0$. 
\end{thm}
\noindent The MDF in the VC theorem allows for many obvious generalization. For the MDF in the Glivenko-Cantelli theorem we give an independent proof based on Hoeffding's inequality. 
\begin{thm}[Glivenko-Cantelli with exponential MDF]\label{thm:Glivenko-Cantelli} 
Given $(X_i)_{i\in \NN}$ i.i.d. in with values in $\RR$ and distribution function $F(x)$ we define 
\[
\hat F_n(x, \om) := \frac{1}{n}\sum_{i=1}^n \ind\{X_i(\om) \lqq x\}, \quad\om \in \Omega, \qquad \mbox{ where }\quad \EE[\hat F_n(x, \cdot)] = F(x),  \quad \mbox{ for all }x \in \RR, ~n\in \NN.
\]
Then 
\[
\lim_{n\ra\infty} \sup_{x\in \RR} |\hat F_n(x, \cdot) - F(x)|= 0\qquad \mbox{a.s. 
with exponential MDF in the sense that:}   
\]
There is $K>0$ such that for any $\e>0$, $\oO_\e = \#\{n\in \NN~|~\sup_{x\in \RR} |\hat F_n(x, \cdot) - F(x)| \gqq \e\}$ and $0 < \eta < \e$ we have 
\[
\EE\big[e^{2\eta^2 \oO_\e}\big]
\lqq K/ (\e^6 (\e -\eta)), \quad \mbox{and in particular }\quad 
\PP(\oO_\e \gqq k)\lqq  K/ (\e^6 (\e -\eta))\cdot \exp(-2 \e^2 k) \quad \mbox{ for all } k\in \NN.
\]
\end{thm}

\begin{proof2}
For $\e\in (0,1)$ set $M = \lceil \frac{1}{\e}\rceil$ 
and $x_0 < x_1 < \dots < x_M$ such that $F(x_i)-F(x_{i+1}) < 1/M$.  
Hence 
\begin{align*}
\PP(\sup_{x\in \RR} |\hat F_n(x, \cdot) - F(x)| \gqq\e) \lqq \sum_{i=1}^M \PP(|\hat F_n(x_i, \cdot) - F(x_i)|\gqq\e).  
\end{align*}
By Hoeffding's inequality \cite{Hoe63}, 
we have $\sup_{x\in \RR}\PP(|\hat F_n(x, \cdot) - F(x)| \gqq \e) \lqq 2 \exp(-2n \e^2)$ for all $\e>0$. Consequently, there is $K>0$ such that for sufficiently small $\e>0$ 
it follows for all $0 < \eta < \e$ that 
\[
\sum_{n=1}^\infty  e^{2\eta^2 n} \cdot \sum_{m=n}^\infty \PP(\hat F_m(x, \cdot) - F(x) \gqq \e)
\lqq (4/\e)e^{-2\e^2} (1- e^{-2 \e^2}) \sum_{n=1}^\infty  e^{2\eta^2n} \cdot  e^{-2 \e^2 n}
\lqq K/(\e^6 (\e-\eta)).
\]
By Corollary~\ref{cor: sufficient poly and exp}, item~(2), 
$\EE[\exp(\e^2 \oO_\e)] \lqq K/(\e^6 (\e-\eta))$ follows 
after an adjustment of $K$. 
 \end{proof2}
 
 \begin{rem}
  The constant $K$ in Theorem~\ref{thm:Glivenko-Cantelli} is optimized in the Dvoretzky-Kiefer-Wolfowitz inequality~\cite{DKW56}, while Sanov's theorem (Theorem~\ref{thm:Sanov} below) optimizes the exponent. 
 \end{rem}

\subsubsection{\textnormal{\textbf{Mean deviation frequency in the Strong Law of Large Numbers}}}We strengthen Etemadi's SLLN 
(see~\cite{Et81} or e.g.\cite{Kle08}, Theorem 15.6) in terms of higher mean deviation frequencies.

    \begin{thm}[Etemadi's SLLN with higher order MDF]
        \label{prop:SLLNhigher}
        Consider an i.i.d.~sequence $(X_i)_{i\in \NN}$ of centered random variables with  
        \begin{equation}\label{e:iid-integrability}\EE[|X_1|^{2q}] <\infty\qquad \mbox{ for some }q\in \NN, q\gqq 1.\end{equation}
        Then we have 
        \[
        \lim_{n\ra\infty } \frac{1}{n} \sum_{i=1}^n X_i  = 0 \qquad \mbox{ a.s. with polynomial MDF of order $0< p < q-1$ in the following sense: }
        \]
        For any $\e>0$ with $\oO_\e := \sum_{n=1}^\infty \ind{\{|\frac{1}{n} \sum_{i=1}^n X_i|\gqq \e\}}$
and $0 < p< q-1$ 
        there is a constant $K_{\e, p}>0$ such that 
        \[
        \EE[\oO_\e^{p}] < p K_{\e, p}, \qquad \mbox{ and in particular }\quad 
        \PP(\oO_\e\gqq k)\lqq p K_{\e,p} \cdot k^{-p}, \qquad \mbox{ for all }k\in \NN.  
        \]
    \end{thm}
    \begin{rem}
(1.) By Lemma~\ref{lem:expected val S_k^2q} below, the maximal order of the overlap statistic $\oO$ is optimal. \\
(2.) The complete independence of $(X_i)_{i\in \NN}$  can be relaxed to the independence of all families $(X_{i_1}, \dots, X_{i_r})$ with indices $I = \{i_1, \hdots, i_r\}\subset \NN$, $2\lqq r\lqq 2q$, which corresponds to the pairwise independence and second moments $(q = 1)$ in Etemadi~\cite{Et81}. 
\end{rem}
\noindent The proof of Theorem~\ref{prop:SLLNhigher} is based on the following lemma.     
 \begin{lem}
    \label{lem:expected val S_k^2q}
    Let $(X_i)_{i=1}^{k_n}$ be a sequence of centered i.i.d. random variables and $S_{k_n}^{2q}:=(\sum_{i=1}^{k_n} X_i)^{2q}$, such that $k_n\gqq q \gqq 1$. If $\Pi$ is the set of all the strictly positive integer partitions
     $\pi:= b_1 + \dots + b_{\ell_{\pi}} =2q$ and $\Pi_{\neq 1}$ the set restricted to $b_m\neq 1$, for all $m=1,\dots , \ell_{\pi}$, then 
    \begin{equation}\label{e:centeredmoments}
    \EE\big[S_{k_n}^{2q}\big]\lqq k_n^q\cdot \frac{(2q)!}{2^{q}}  \sum_{\pi\in \Pi_{\neq 1}} \prod_{m=1}^{\ell_{\pi}} \EE[X^{b_m}_1].\end{equation}
    In fact, the order $k_n^q$ in \eqref{e:centeredmoments} is optimal, whenever $X_1\neq 0$ a.s.  
\end{lem}

\begin{proof2}
    Since the random variables are independent, by applying $\EE$ to multinomial expansion of $S_{k_n}^{2q}$ we get 
    \begin{align}
        \EE\big[S_{k_n}^{2q}\big] = \sum_{b_1 + \dots + b_{k_n}=2q} \binom{2q}{b_1, \dots , b_{k_n}} \prod_{m=1}^{k_n} \EE[X_{i_m}^{b_m}].
        \label{eq:redundant notation} 
    \end{align}
Note that the random variables are centered and identically distributed. Hence, the left-hand side of \eqref{eq:redundant notation} only depends on $b_m\notin \{0, 1\}$, $m\in \{1, \dots, k_n\}$ and the subscripts in $\EE[X^{b_m}_{i_m}]$ can be set equal to $1$. 
    This leaves us with counting how many times a partition $\pi \in \Pi_{\neq 1}$ apprears in the exponents of $\prod_{m=1}^{k_n} \EE[X^{b_m}_1]$. If $\pi$ has $\ell_{\pi}$ terms, then it will appear $k_n!/(k_n-\ell_{\pi}-1)!$ times, as we just need to ensure that the $\ell_{\pi}$ subindices are different. Therefore, \eqref{eq:redundant notation} is bounded by  
    $$\EE\big[S_{k_n}^{2q}\big] \lqq \sum_{\pi \in \Pi_{\neq 1}} \binom{2q}{b_1, \dots , b_{\ell_\pi}} k_n ^{\ell_\pi} \prod_{m=1}^{\ell_\pi} \EE[X^{b_m}_1]=: \sum_{\pi \in \Pi_{\neq 1}} c_{\pi} \prod_{m=1}^{\ell_\pi} \EE[X^{b_m}_1].$$
    \noindent Furthermore, by the pigeonhole principle, the partitions in $\Pi_{\neq 1}$ can have at most $q$ terms without any $b_m\in \{0, 1\}$. In particular, there is just one partition $\pi^*$ such that $\pi^*\in \Pi_{\neq 1}$ and $\ell_{\pi^*} =q$, which is that with $b_m=2$ for $m=1, \dots , q$. We now show that for $k_n\gqq 2$, $b_{\pi^*} = \max\{b_{\pi}~|~\pi \in \Pi_{\neq 1}\}$. To do so, consider any partition $\pi\in \Pi_{\neq 1}$ of length $\ell_\pi<q$. We have that $b_m\gqq 2$ for $m=1, \dots , \ell_\pi$, and therefore
    \begin{align*}
        \binom{2q}{b_1, \dots , b_{\ell_\pi}} = \frac{(2q)!}{b_1! \dots b_{\ell_\pi} !} \lqq \frac{(2q)!}{2^{\ell_\pi}},
    \end{align*}
        which implies \eqref{e:centeredmoments}, since for any $n\in \NN$ we have the equivalences 
    \begin{align*}
        \frac{(2q)!}{2^{\ell_\pi}}\cdot k_n ^{\ell_\pi} \lqq \frac{(2q)!}{2^{q}} \cdot k_n ^{q} = \binom{2q}{2, \dots , 2}\cdot k_n ^{q} \quad \Leftrightarrow \quad 2^{q-\ell_\pi} \lqq k_n ^ {q-\ell_\pi} \quad \Leftrightarrow \quad 2 \lqq k_n.\\[-1cm]
    \end{align*} 
\end{proof2}
   
    \begin{proof2} (of Theorem~\ref{prop:SLLNhigher})
        The fact that $(X_i)_{i\in \NN}$ fulfills the SLLN follows directly from Theorem 5.16 in~\cite{Kle08}. 
For $S_n:= \sum_{i=1}^n X_i$, $n\in \NN$, the proof there boils down to an application of Markov's inequality and the subsequence $(k_n)_{n\in \NN}$ given by $k_n:= \lfloor (1+\e)^n \rfloor\gqq \frac{1}{2}(1+\e)^n$. In particular, for $\psi(x)=x^{2q}$ we obtain by Markov's inequality that
    \begin{align}
        \PP\Big( \Big|\frac{S_{k_n}}{k_n}\Big| \gqq \e  \Big)
        \lqq \psi\Big(\Big|\frac{S_{k_n}}{k_n} \Big|\Big) \cdot\psi^{-1}(\e)
        =\Big(\frac{1}{\e \cdot k_n}\Big)^{2q} \EE[S_{k_n}^{2q}]. 
        \label{e:Markov ineq S_kn^2q}
    \end{align}
    Lemma \ref{lem:expected val S_k^2q} yields a constant $K_{q}>0$ 
    such that $\EE[S_{k_n}^{2q}]\lqq K_q \cdot k_n^{q}$ for all $n\in \NN$. 
    Hence, \eqref{e:Markov ineq S_kn^2q} is bounded by
    \begin{align}
        \PP\Big( \Big|\frac{S_{k_n}}{k_n}\Big| \gqq \e  \Big) \lqq \frac{K_q}{(\e ^{2}\cdot k_n)^q}, \qquad n\in \NN 
\label{e: bound S_kn/k_n}
    \end{align}
We now pass from the subsequence $(k_n)_{n\in \NN}$ to the general index $\ell\in \NN$. We start by bounding 
\begin{align}\label{e: piv1}
\Big| \frac{S_\ell}{\ell} \Big| &\lqq  \Big| \frac{S_{k_n}}{k_n} \Big| +  \Big| \frac{S_\ell}{\ell}- \frac{S_{k_n}}{k_n} \Big| \qquad \mbox{ for }
k_n\lqq \ell \lqq k_{n+1}.  
\end{align}    
Moreover, for a sufficiently large $n\in \NN$, we have that $k_{n+1}\lqq (1+2\epsilon) k_{n}$ and for every $\ell \in \{ k_n, \dots , k_{n+1}\}$ 
        \begin{align*}
          &  \frac{1}{1+2\epsilon}\frac{S_{k_n}}{k_n}\lqq \frac{S_{k_n}}{k_{n+1}}\lqq \frac{S_{\ell}}{\ell} \lqq \frac{S_{k_{n+1}}}{k_n}\lqq (1+2\e) \frac{S_{k_{n+1}}}{k_{n+1}} ~\vzv~ \frac{-2\e}{1+2\epsilon}\frac{S_{k_n}}{k_n}\lqq  \frac{S_{\ell}}{\ell} - \frac{S_{k_n}}{k_n}\lqq  (1+2\e) \frac{S_{k_{n+1}}}{k_{n+1}}-\frac{S_{k_n}}{k_n}. 
\end{align*}
Consequently
         \begin{align}
         \max_{k_n\lqq \ell \lqq k_{n+1}} &\Big| \frac{S_\ell}{\ell}- \frac{S_{k_n}}{k_n} \Big|  
          \lqq \max\{\Big(\frac{2\e}{1+2\epsilon}\Big)|\frac{S_{k_n}}{k_n}|, |(1+2\e) \frac{S_{k_{n+1}}}{k_{n+1}}-\frac{S_{k_n}}{k_n}| \}\nonumber \\
          &\lqq \Big(1-\frac{1}{1+2\epsilon}\Big)|\frac{S_{k_n}}{k_n}| + |(1+2\e) \frac{S_{k_{n+1}}}{k_{n+1}}-\frac{S_{k_n}}{k_n}|\lqq 2\e|\frac{S_{k_n}}{k_n}| + |(1+2\e) \frac{S_{k_{n+1}}}{k_{n+1}}-\frac{S_{k_n}}{k_n}|.\label{e: piv2}
          \end{align} 
Note that $k_{n+1} \lqq (1+2\e) k_n$, such that $k_{n+1} - k_n \lqq 2\e k_n$. 
Combining this with \eqref{e: piv1} and \eqref{e: piv2} we obtain          
\begin{align}
\sum_{\ell=1}^\infty \ell^{p} \sum_{m=\ell}^\infty \PP\Big( \Big|\frac{S_{m}}{m} \Big| \gqq \e  \Big) &\lqq \sum_{n=1}^\infty \sum_{\ell =k_{n}}^{k_{n+1}-1} \ell^{p}\sum_{m=k_n}^\infty \PP\Big( \Big|\frac{S_{m}}{m} \Big| \gqq \e  \Big)\nonumber\\
 &\hspace{-2cm}\lqq \sum_{n=1}^\infty k_{n+1}^{p} (k_{n+1}- k_{n}) \sum_{m=k_n}^\infty \PP\Big( \Big|\frac{S_{m}}{m} \Big| \gqq \e  \Big) \lqq 2\e (1+2\e)^p \sum_{n=1}^\infty k_{n}^{p+1} \sum_{m=k_n}^\infty \PP\Big( \Big|\frac{S_{m}}{m} \Big| \gqq \e  \Big).\label{e: summabilitycrit}
 \end{align}
On the other hand, thanks to \eqref{e: piv2} we get the upper bound 
\begin{align}
\sum_{m=k_n}^\infty \PP\Big( \Big|\frac{S_{m}}{m} \Big| \gqq \e  \Big)
&\lqq \sum_{t=n}^\infty \sum_{s=1}^{k_{t+1}-k_{t}} 
\PP\Big( \Big|\frac{S_{k_t+s}}{k_t+s} \Big| \gqq \e  \Big)\nonumber\\
&\lqq 2 \e \sum_{t=n}^\infty k_t  \Big(\PP(|\frac{S_{k_t}}{k_t}|\gqq\frac{\e}{3}) + \PP(2\e|\frac{S_{k_t}}{k_t}|\gqq\frac{\e}{3}) 
+ \PP(|(1+2\e) \frac{S_{k_{t+1}}}{k_{t+1}} - \frac{S_{k_t}}{k_t}|\gqq \frac{\e}{3})
\Big).\label{e: upperbound}
\end{align}
Furthermore, for $0 <\e <1/6$ and $0 < \ti \e < \e/8 \lqq (\e / 6) (1+2\e)^{-1}$ we apply \eqref{e: bound S_kn/k_n} to \eqref{e: upperbound}and obtain that 
\begin{align}\label{e: finalsum}
\sum_{m=k_n}^\infty \PP\Big( \Big|\frac{S_{m}}{m} \Big| \gqq \e  \Big)
\lqq 8 \e \sum_{t=n}^\infty k_t \PP\Big(|\frac{S_{k_t}}{k_t}\Big|\gqq\ti \e\Big) 
\lqq \frac{8^{q+1} K_q}{\e^{2q-1}} \sum_{t=n}^\infty \frac{1}{k_t^{q-1}} 
\lqq \frac{8^{q+1} }{\e^{2q-1}}  \frac{K_q(1+\e)^{q-1}}{(1+\e)^{q-1}-1} \frac{1}{k_n^{q-1}}
=: \frac{\hat K_{q,\e}}{k_n^{q-1}}. 
\end{align}
Comining \eqref{e: summabilitycrit} with \eqref{e: finalsum} we obtain 
\begin{align*}
\sum_{\ell=1}^\infty \ell^{p} \sum_{m=\ell}^\infty \PP\Big( \Big|\frac{S_{m}}{m} \Big| \gqq \e  \Big) 
\lqq 2\e (1+2\e)^p \hat K_{q, \e} \cdot \sum_{n=1}^\infty k_n^{p+2-q} <\infty 
\end{align*}
for all $p < q -2$ the right-hand side is finite and Corollary~\ref{cor: aspoloutlier} yields for $\oO_\e := \sum_{n=1}^\infty \ind{\{|S_{n}/n| \gqq \e\}}$
the desired estimate 
$\EE[\oO_\e^{p+1}] < \infty$. Renaming the constants finishes the proof. 
\end{proof2}

\subsubsection{\textnormal{\textbf{MDF under large deviation principles, Cram\'er and Sanov's theorem}}}\hfill

\begin{thm}[A large deviations principle implies exponential MDF]\label{exp:ldp}
Assume a family $(\mu_n)_{n\in \NN}$ of probability measures with $\mu_n(A) = \PP(E_n(A))$, $n\in\NN$, for some $E_n(A)\in \bB(\RR^d)$ on a common probablity space $(\Omega, \aA, \PP)$ which satisfies the following upper bound of a large deviations principle (LDP) for some good rate function $\jJ$:  
\begin{align}\label{e: LDP} 
\limsup_{n\ra\infty} \frac{1}{n} \ln \PP(E_n(A)) \lqq - \inf_{x\in \bar A} \jJ(x), \qquad A\in \aA.  
\end{align}
For $A\in \aA$ consider the overlap statistic $\oO_A$ of $(E_n(A))_{n\in \NN}$. 
Then for all $A\in \aA$ there is $C_A>0$ such that 
for any $0 < p < \inf_{x\in \bar A} \jJ(x)$ we have 
\[
\EE[e^{p\oO_A}] 
\lqq C_A \Big(1-\exp(-(\inf_{x\in \bar A} \jJ(x)))\Big)^{-1}\Big(1-\exp(-(\inf_{x\in \bar A} \jJ(x)-p)\Big)^{-1}  = K_A < \infty.
\]
In particular, we have for any $0 < p < \inf_{x\in \bar A} \jJ(x)$ we have that  
$\PP(\oO_A\gqq k) \lqq K_A\exp(-pk)$ for all $k\in \NN$. 
\end{thm}

\begin{proof2} Consider $A\in \aA$ and $0 < p < \inf_{x\in \bar A} \jJ(x)$. Combining \eqref{e:expmoment} and \eqref{e: LDP} yields $C_A>0$ such that 
$\mu_n(A) = \PP(E_n(A)) \lqq C_A \exp(-n \inf_{x\in \bar A} \jJ(x))$ for all $n\in \NN$. 
The result follows from Corollary~\ref{cor: sufficient poly and exp}, item~(2). 
\end{proof2}
\noindent The following example is motivated in comparative analysis of DNA sequence matching (see \cite{AGW90} and \cite{DZ98}, p.83).  
\begin{exm}[Exponential MDF for long rare segments in random walks, \cite{DZ98} Sec.~3.2]\label{exm:longraresegments}  
Consider a random walk $(S_n)_{n\in \NN}$, $S_0 = 0$ and $S_n = \sum_{i=1}^n X_i$, $n\in \NN$ for an i.i.d. sequence $(X_i)_{i\in \NN}$ with values in $\RR^d$. For a Borel set $B\subset \RR^d$ we consider the maximal segment length $R_n$ of $(S_n)_{n\in \NN}$ whose empirical mean belongs to $B$, and the first occurence time $\tau_r$ of a such a segment of length $r$, respectively as  
\begin{align*}
R_n := \max\{\ell-k~|~0\lqq k\lqq \ell\lqq n: \frac{S_\ell-S_k}{\ell - k}\in B\},
\quad \tau_r := \inf\{\ell\in \NN~|~ \frac{S_\ell-S_k}{\ell - k}\in B \mbox{ for some }0\lqq k\lqq \ell -r\}. 
\end{align*}
We assume the following strong version of a large deviation principle 
\begin{align}\label{e:strongLDP}
\jJ(B) := \lim_{n\ra\infty}(1/n) \ln \PP(S_n/n\in B).
\end{align}
Then the following a.s. limits are both valid with exponential MDF: 
\begin{align}\label{e: long rare}
\lim_{r\ra\infty} \frac{\ln \tau_r}{r}= \jJ(B) ~\mbox{ a.s.} \qquad \mbox{ and }\qquad \lim_{n\ra\infty} \frac{R_n}{\ln(n)} = \frac{1}{\jJ(B)} ~\mbox{ a.s.}
\end{align}
More precisely, if we define the one-sided deviation frequencies for $\e>0$ 
\begin{align*}
\oO_\e^+ &:= \sum_{n=1}^\infty\ind\{R_n/ \ln(n)  \gqq (\jJ(B) -\e)^{-1}\}, \qquad \oO_\e^- := \sum_{n=1}^\infty\{ R_n/ \ln(n)  \lqq (\jJ(B) +\e)^{-1}\}, \\
\uU_\e^+ &:= \sum_{r=1}^\infty \ind\{\ln(\tau_r) / r  \gqq \jJ(B) +\e\}, \hspace{1,7cm} 
\uU_\e^- := \sum_{r=1}^\infty \{\ln(\tau_r)/r  \lqq \jJ(B) -\e\}.
\end{align*}
Then for all $0 <\eta <\e$ there are constants $\kK_1(\e, \eta)>0$, $p>0$ and $\kK_2(\e, p)>0$ such that 
\begin{align*}
\EE\big[e^{\eta \oO_\e^+}\big] = \EE\big[e^{\eta \uU_\e^-}\big]\lqq \kK_1(\e, \eta)+1 \qquad \mbox{ and }\qquad \EE\big[e^{p e^{\eta \oO_\e^-}}\big] = \EE\big[e^{p e^{\eta \uU_\e^+}}] \lqq \kK_2(\e, p).
\end{align*}
In particular, the Markov inequality yields for all $k\in \NN$ the exponential and doubly    exponential bounds 
\begin{align*}
\qquad &~\PP(\oO_\e^+\gqq k) = \PP(\uU_\e^-\gqq k)\lqq 
\kK_1(\e, \eta) e^{-\eta k}\quad\qquad  \mbox{ and }\qquad\, \PP(\oO_\e^-\gqq k) =\PP(\uU_\e^+\gqq k)\lqq 
\kK_2(\e, \eta) e^{- p e^{\e k}}. 
\end{align*}
\textbf{Sketch of proof: }We discuss the case $0 < \jJ(B) <\infty$, and $\jJ(B)= \infty$ can be treated similarly. We show the result for~$\tau_r$. The complementary results for $R_m$ are then a consequence of the duality $\{R_n\gqq r\} = \{\tau_r\lqq n\}$ for all $r, n\in \NN$. In the proof of Theorem~3.2.1 on p. 84 in \cite{DZ98} the authors apply  \eqref{e:strongLDP} and show the estimate 
\begin{align*}
\PP(\tau_r \lqq m) \lqq m \sum_{n=r}^\infty \PP(S_n / n \in B).\\[-7mm]  
\end{align*}
Then, for a fixed $\e>0$ there are positive constants $c = c(\e), \ti c = \ti c(\e)>0$ 
such that for $m = \lfloor \exp(r (\jJ(B)-2\e))\rfloor$ 
\begin{align*}
 \sum_{r=\varrho}^\infty  \PP(\tau_r \lqq \exp(r (\jJ(B) -\e))) 
\lqq \sum_{r=\varrho}^\infty  \exp(r (\jJ(B)-2\e)) \sum_{n=r}^\infty c \exp(-n(\jJ(B)-\e))
\lqq (2\ti c /\e) \exp(-\varrho \e).  
\end{align*}
Now, for any $0 < \eta <\e$ we have 
\begin{align*}
\kK_1(\e, \eta) :=\sum_{\varrho = 1}^\infty  \exp(\varrho \eta) \sum_{r=\varrho}^\infty  \PP(\tau_r \lqq \exp(r (\jJ(B) -\e))) \lqq (2\ti c/\e) \sum_{\varrho = 1}^\infty  \exp(-\varrho (\e-\eta)) <\infty. 
\end{align*}
Hence by Corollary~\ref{cor: sufficient poly and exp}, item~(2), we have 
$\EE[e^{\eta \uU_\e^-}]\lqq \kK_1(\e, \eta) +1< \infty$. For the upper bound an application of \eqref{e:strongLDP} given on p.85 of \cite{DZ98} yields that $\PP(\tau_r \gqq m)\lqq \exp(-\lfloor m/r\rfloor \PP( S_n /n\in B))$ for all $r, m \in \NN$. Hence for $m = \lfloor \exp(r (\jJ(B)+2\e))\rfloor$ there are constants $c_1, c_2, c_3>0$ such that 
\begin{align*}
\sum_{r = \varrho}^\infty \PP(\tau_r > \exp(r(\jJ(B) +\e)))\lqq \sum_{r = \varrho}^\infty \exp(- (c_1/r) e^{\e r}) \lqq c_3\sum_{r = \varrho}^\infty \exp(- c_2 e^{\e r}) 
\end{align*}    
and for all $0 <p < c_2$ we have
\begin{align*}
\kK_2(\e, p) :=\sum_{\varrho = 1}^\infty \exp( p e^{\e \varrho}) 
\sum_{r = \varrho}^\infty \PP(\tau_r > \exp(r(\jJ(B) +\e)))
&\lqq \sum_{\varrho = 1}^\infty \exp(-(c_2- p) e^{\e \varrho})  < \infty. 
\end{align*}
Hence for $\sS(N) := \sum_{n=1}^N \exp( p e^{\e n})$ 
Theorem~\ref{prop: sufficient} implies that $\EE[\exp( p e^{\e \uU_\e^+})] \lqq \EE[\sS(\uU_\e^+)] \lqq \kK_2(\e, p) <\infty$.   
\end{exm}

\noindent In the sequel we apply the SLLN in the formulation of Cram\'er and Sanov's theorem. 

\begin{thm}[Cram\'er's theorem with exponential MDF]\label{thm: Cramer}
Consider an i.i.d. family $(X_i)_{n\in \NN}$ with values in $\RR^d$. Then  
$\Lambda(\la) := \ln( \EE[e^{\lgl \la, X_1\rgl}]) < \infty$ for all $\la\in \RR^d$ implies that 
\[
\lim_{n\ra\infty} \frac{1}{n} \sum_{i=1}^n X_i =  \EE[X_1] \quad \mbox{a.s. with exponential MDF in the following sense:} 
\]
For any $\e>0$, $\oO_\e := \sum_{n=1}^\infty \ind{\{|\frac{1}{n} \sum_{i=1}^n X_i - \EE[X_1]| \gqq \e\}}$ and $0 < p < \inf_{|\la^*|>\e} \Lambda^*(\la^*)$ there is a constant $C_{\e, p}>0$ such that 
\begin{align*}
\EE\big[e^{p \oO_\e}\big] \lqq C_{\e, p}
\Big(1-\exp\big(-\inf_{|\la^*|>\e} \Lambda^*(\la^*)\big)\Big)^{-1}
\Big(1-\exp\big(-(\inf_{|\la^*|>\e} \Lambda^*(\la^*)-p\big)\Big)^{-1}  =: K_{\e, p} < \infty,
\end{align*}
where the $\La^*$ is the Fenchel-Legendre transform of $\La$. Additionally  
$\PP(\oO_\e\gqq k) \lqq K_{\e, p}\cdot \exp\big(-p k\big)$ for all $k\in \NN.$ 
\end{thm}
\noindent The proof is a direct application of Example~\ref{exm: exp}.\\ 
\paragraph{\textbf{The setting of Sanov's theorem: }} Consider an i.i.d. sequence $X = (X_i)_{i\in \NN}$ with values in a Polish space $\Sigma$ and a common distribution $\mu$. 
We equip $ \mbox{Prob}(\Sigma)$ with the so-called $\tau$-topology  (see \cite{DZ98}, Sec.~6.2, p.~263 or \cite{Ce07}, Chapter 23). 
We denote the empirical law of $X$ by\\[-5mm]
\begin{align*}
\lL^X_n := \frac{1}{n} \sum_{i=1}^n \delta_{X_i}  \in \mbox{Prob}(\Sigma).\\[-7mm]
\end{align*}
and the \textit{relative entropy}, also known as \textit{Kullback-Leibler divergence} between $\mu, \nu \in \mbox{Prob}(\Sigma)$, as  
\begin{align*}
D_{\mathrm{KL}}(\nu |\mu) := \int_{\Sigma} f \ln(f) d\mu,  \mbox{ if the Radon-Nikodym derivative } f = \frac{d\nu}{d\mu}\mbox{ exists,  and equal to } \infty \mbox{ otherwise.}  
\end{align*}

\begin{thm}[Sanov's theorem with exponential MDF]\label{thm:Sanov} 
Under the preceding setting the empirical measures satisfy for any $\tau$-measurable set $B\subset \mathrm{Prob}(\Sigma)$ the upper bound of the LDP 
\[
\limsup_{n\ra\infty} \frac{1}{n} \ln \PP(\lL^X_n \in B)\lqq \inf_{\nu \in \bar B} D_{\mathrm{KL}}(\nu|\mu).
\]
Then for $E_n(B) := \{\lL^X_n \notin B\}$, $\oO_B := \sum_{n=1}^\infty \ind_{E_n(B)}$ and 
$0 < p < \inf_{\nu \in \bar B} D_{\mathrm{KL}}(\nu|\mu)$ there is $C_{p, B}>0$ such that  
\begin{align*}
\EE\big[e^{p\oO_B}\Big] 
&\lqq C_{p, B}\cdot \Big(1- \exp\big(-\inf_{\nu\in \bar B} D_{\mathrm{KL}}(\nu|\mu)\big)\Big)^{-1}  
  \cdot \Big(1-\exp\big(-(\inf_{\nu\in \bar B} D_{\mathrm{KL}}(\nu|\mu)-p)\big)\Big)^{-1} = K_{p, B}, 
\end{align*}
and $\PP(\oO_B\gqq k) \lqq K_{p, B}\cdot \exp\big(-p k \big)$ for all $k\in \NN_0$. 
The proof is a direct application of Theorem~\ref{exp:ldp}.

\end{thm}
\noindent 

\subsubsection{\textnormal{\textbf{MDF quantification of the Method of Moments}}}\hfill

\begin{thm}[Method of Moments]\label{thm:MDFinMOM}
Consider an i.i.d. sequence $(X_i)_{i\in \NN}$ with values in $\RR$ and common distribution $\PP_\theta$ that depends on an unknown vector of parameters, $\theta = (\theta_1, \dots, \theta_k) \in \Theta$, for some open subset $\Theta \subset \RR^k$, $k\in \NN$. 
We set $\theta \mapsto M(\theta) := (m^1(\theta), \dots, m^k(\theta))$, where $m^j(\theta) := \EE[X_1^j(\theta)], 1\lqq j\lqq k$, and define 
\begin{align*}
&\bar X_n(\theta) := (\bar X_{n, 1}(\theta), \dots, \bar X_{n,k}(\theta)), \qquad \mbox{ with }\quad 
\bar X_{n,j}(\theta) :=  \frac{1}{n} \sum_{i=1}^n X^j_{i}(\theta).\\[-10mm]
\end{align*}
\textbf{Hypotheses:}
\begin{enumerate}\item[\textbf{(i)}] Let $\sup_{\theta\in \Theta}\EE[|X_1(\theta)|^{2q \cdot k}] < \infty$ for some $q\in\NN$, $q\gqq 2$. 
 \item[\textbf{(ii)}] For any $\theta_0\in \Theta$ the mapping $\Theta \ni \theta \mapsto M(\theta)$ is continuous, bijective, and its inverse $M^{-1}$ is continuously differentiable in an open neighborhood of $\theta_0\in \Theta$.
\item[\textbf{(iii)}] $M(\theta)$ only depends on the odd powers $m^j$, $j=2i -1$, $1\lqq i\lqq \lfloor k/2\rfloor$. 
\end{enumerate}
Then for any $\theta_0\in \Theta$ the estimator
$\hat \theta_n := M^{-1}(\bar X_n(\theta_0))$, $n\in\NN$, 
satisfies $\hat \theta_n \ra \theta_0$ a.s. as $n\ra\infty$ with $p$-th MDF for $0< p< q-1$. For any such $p$ and $\e>0$ there is $K_{\e, p}>0$ such that 
$\oO_\e:= \sum_{n=1}^\infty \ind{\{|\hat \theta_n- \theta_0|\gqq\e\}}$ satisfies 
\begin{align*}
\EE[\oO_\e^{p}] \lqq p K_{\e, p}\qquad\mbox{ and } \qquad 
\PP(\oO_\e \gqq k)\lqq p K_{\e, p} \cdot k^{-p}
\qquad \mbox{ for all }k\in \NN. 
\end{align*}
\end{thm}
\noindent Note that exponential integrability on the $X_i$ the result can be improved according to Theorem~\ref{thm: Cramer}.
\begin{proof2} 
Fix $\theta_0\in \Theta$. Then for $\e>0$ sufficently small we have $\la = \la(\e) = \mbox{specrad}(D_{\theta_0} M)-\e>0$ such that 
 \begin{align*}
 E_n(\e) &:= \{|\hat \theta_n- \theta_0|\gqq \e\} 
 = \{M^{-1}(\bar X_n(\theta_0)) \notin B_{\e}(\theta_0)\}= \{\bar X_n(\theta_0) \notin M(B_{\e}(\theta_0))\}\\[2mm]
 &\subset  \{\bar X_n(\theta_0) \in (D_{\theta_0} M) B_{\e/2}^c(M(\theta_0)))\} 
 \subset \{\bar X_n(\theta_0) \in B_{\la\cdot \e/2}^c(M(\theta_0)))\}= \{|\bar X_n(\theta_0)-M(\theta_0)|\gqq \la \cdot \e/2\}. 
 \end{align*}
 \noindent We denote $F_n(\e):= \{|\bar X_n(\theta_0)-M(\theta_0)|\gqq \la\cdot \e/2\}$, which results in  
 $\oO_\e = \sum_{n=1}^\infty \ind(E_n(\e)) \lqq \sum_{n=1}^\infty \ind(F_n(\e))=:\uU_\e$ by monotonicity. Theorem~\ref{prop:SLLNhigher} then implies for $0 < p < q-1$ that 
 $\EE[ \oO_\e^{p}]\lqq \EE[\uU_\e^{p}] <\infty$.
 \end{proof2}

\medskip 
\subsection{MDF in the Law of the Iterated Logarithm and strong explicit schemes of SDE}\hfill\\[-6mm]
\subsubsection{\textnormal{\textbf{MDF in the Law of the Iterated Logarithm for Brownian motion}}}The law of the iterated logarithm provides a.s. asymptotic bounds for the growth of the trajectories of a standard Brownian motion. It is natural to ask, how many times this asymptotically strong bound will be trespassed along a given sequence $\al^n$ of partitions, until the a.s. asymptotics kicks in. 

\begin{thm}[Exceedance frequency in the law of the iterated logarithm]\label{thm: lil} 
Let $W = (W_t)_{t\gqq 0 }$ be a standard Brownian motion with values in $\RR$. Then we have 
\[
\limsup_{t\ra \infty} \frac{W_t}{\sqrt{2t \ln(\ln(t))}} \lqq 1 \qquad \mbox{a.s.}
\]
Moreover, for any $\al>1$ there is a constant $\kK_\al>0$ such that the exceedance frequency, given by  
\[
\oO_\al = \sum_{n=1}^\infty \ind_{E_n(\al)}, \qquad E_n(\al):= \{\sup_{t\in (\al^n, \al^{n+1}]} W_t > \sqrt{\al} \sqrt{2 \cdot \al^n \ln(\ln(\al^n))}\},\quad n\in \NN, 
\]
is finite a.s. and by Example~\ref{ex: poly} here exists a constant $\kK_\al>0$ such that 
$\EE[\oO_\al^{1+(\al-2)\vee 0}] \lqq \kK_\al$ and 
\[
\PP(\oO_\al\gqq k)\lqq
\kK_\al k^{-(1+(\al-2)\vee 0)}
\qquad \mbox{ for all }k\in\NN, \quad k\gqq 1.
\]
\end{thm}

\begin{proof2}
Set $\al>1$ and $t_n = \al^n$ and $f_\al(t) := 2 \al^2 \ln(\ln(t))$,  so we have $\sqrt{t_n f_\al(t_n)} = \sqrt{2 \al t_{n+1} \ln(\ln(t_n))}$. Formula (22.3) on p. 496 of \cite{Kle08} 
states the existence of a constant $C>0$ such that for all $n\in\NN$ 
\begin{align*}
\PP(E_n)
\lqq  C n^{-\al} \quad \mbox{ and thus }
\sum_{n=1}^\infty \PP(E_n(\al))
\lqq C \zeta(\al)
=: c_\al < \infty.
\end{align*}
Consequently, the first Borel-Cantelli Lemma given by \eqref{first BC} yields that for all $\al>1$ we have 
\[
\limsup_{n\ra\infty} \sup_{t\in [t_n, t_{n+1}]} \frac{W_t}{\sqrt{2 \al t \ln(\ln(t))}}\lqq 1, \qquad \mbox{ hence }
\qquad \limsup_{t\ra\infty} \frac{W_t}{\sqrt{2 \al t \ln(\ln(t))}}\lqq 1 \mbox{ and by \eqref{e:expectation} }\quad \EE[\oO_\al] \lqq c_\al. 
\]
Moreover, by Example \ref{ex: poly} there exists $\kK_\al>0$ such that $\EE[\oO_\al^{1+(\al-2)\vee 0}] \lqq \kK_\al$.
\end{proof2}
\begin{rem}
\begin{enumerate}
 \item Analogous bounds can be derived for LIL
 of random walks~\cite{Kle08} and $\alpha$-stable processes~\cite{Bertoin}. 
 \item With more technical effort, Theorem~\ref{thm: lil} can be generalized to other diverging sequences $(t_n)_{n\in \NN}$.
\end{enumerate}
\end{rem}
\subsubsection{\textnormal{\textbf{A MDF error estimate of strong explicit numerical schemes for SDE}}} 
On a given probability space $(\Omega, \aA, \PP)$ and a time interval $[0, T]$, $T>0$, we consider the scalar stochastic differential equation 
\begin{equation}\label{e: SDE}
\mathrm{d}X = a(t, X) \mathrm{d}t + b(t, X) \mathrm{d}W, \qquad X_0 = Z_0
\end{equation}
for a scalar standard Brownian motion $(W_t)_{t\gqq 0}$. 
Under standard Lipschitz and boundedness conditions on the coefficientes $a$ and $b$ (see for e.g. \cite{KS98}) and square inegrability of $Z_0$, equation \eqref{e: SDE} has a unique strong solution. In \cite{KP92}, Sec.~11.2, the authors derive the following explicit $1.5$-order strong scheme:  
Fix $Y_0 = Z_0$. For $N\in \NN$ we consider a partition $0 = \tau_0 < \tau_1 < \dots < \tau_N = T$ with $\delta := \sup_{n}(\tau_{n+1}-\tau_n)$ and the piecewise linear approximation $(Y_t)_{t\in [0, T]}$ of $(X_{t})_{t\in [0, T]}$. In particular, for $n_t := \max\{n\in \{0, \dots, N\}~|~\tau_n\lqq t\}$, $t\in [0, T]$, and 
\begin{align*}
&Y_t = Y_{n_t} + \frac{t-\tau_n}{\tau_{n_t+1}- \tau_{n_t}}(Y_{n_t+1} - Y_{n_t}), \qquad \mbox{ we define the recursive scheme for }(Y_n)_{n=0, \dots, N} \mbox{ given by}\nonumber\\
&Y_{n+1} = Y_n + b \Delta W + \frac{\Delta Z}{2\sqrt{\Delta}}\Big(a(\Upsilon_+)-(\Upsilon_-)\Big)  + \frac{\Delta}{4} \Big(a(\Upsilon_+)-a(\Upsilon_-)\Big)  + \frac{(\Delta W)^2-\Delta)}{4\sqrt{\Delta}} \Big(b(\Upsilon_+)-b(\Upsilon_-)\Big) \\
&+ \frac{(\Delta W)^2-\Delta}{2\Delta} \Big(b(\Upsilon_+)-b(\Upsilon_-)\Big) \big(\Delta W \cdot \Delta - \Delta Z\big)+ \frac{\Delta W}{4\sqrt{\Delta}} \Big(b(\Phi_+)-b(\Phi_-) - b(\Upsilon_+) + b(\Upsilon_-)\Big) \big(\frac{1}{3} (\Delta W)^2 -\Delta\Big),
\end{align*}
and the notation $\Delta = \tau_{n+1}-\tau_n$, $\Delta W = W_{\tau_{n+1}} - W_{\tau_n}$, $\Upsilon_\pm = Y_n + a \Delta \pm b \sqrt{\Delta}$ 
and $\Phi_\pm = \Upsilon_+ \pm b(\Upsilon_+)  \sqrt{\Delta}$. 
Furthermore, if $\delta = \delta_N \lqq C T/N$, $N\in \NN$, for some $C>0$, Theorem~11.5.1 in \cite{KP92} yields $\mathcal{K}_T>0$~such~that   
\begin{equation}\label{e: summable scheme}
\EE[|X(T)- Y^{\delta_N}_T|] \lqq \mathcal{K}_T \delta_N^{3/2} = \mathcal{K}_T (CT)^{3/2} \cdot N^{-3/2}, \qquad N\in \NN. 
\end{equation}
The summability of the right-hand side of \eqref{e: summable scheme} in $N$ and \eqref{e:expectation} yield for $\e>0$ and 
$\oO_\e := \sum_{n=1}^\infty \ind{\{|X(T)-  Y^{\delta_N}_T|\gqq\e\}}$ that $\EE[\oO_\e] \lqq K_1$,  
and hence 
$\PP(\oO_\e \gqq k)\lqq K_1 \cdot k^{-1}, k\in \NN$.
The rate of convergence is better 
for (more involved) higher order schemes, see Theorem~11.5.2 in \cite{KP92}. 
For example, the scheme of strong order $\gamma = 5/2$ given in Corollary~\ref{cor: aspoloutlier} of~\cite{KP92} enjoys the algebraic MDF error estimate 
\begin{align*}
\PP(\#\{N\in \NN~|~|X(T)-  Y^{\delta_N}_T|\gqq \e\}\gqq k)\lqq (5/2) K_{3/2} \cdot k^{-3/2}, \qquad k\in \NN.\\[-7mm]
\end{align*}

\section*{Acknowledgments} 

\noindent MAH acknowledges support by project INV-2019-84-1837 of Facultad de Ciencias at Universidad de los Andes.


\begin{thebibliography}{99}
\bibitem{Fr73}
Freedman,~D.,1973. 
\newblock Another note on the Borel-Cantelli lemma and the strong law with the Poisson approximation as a by-product
\newblock Ann. Probab., Vol. 1, No. 6, 910--925.
\end{thebibliography}
\end{document}